\documentclass[12pt,twoside,leqno]{article}
\usepackage{latexsym}
\def\implies{\; \Longrightarrow \;}

\begin{document}
\begin{center}
{\bf ON THE SEMIGROUP OF PARTIAL ISOMETRIES OF A FINITE CHAIN}\\[4mm]
\textbf{R. Kehinde and A. Umar}\\

\end{center}

\newtheorem{theorem}{{\bf Theorem}}[section]
\newtheorem{prop}[theorem]{{\bf Proposition}}
\newtheorem{lemma}[theorem]{{\bf Lemma}}
\newtheorem{corollary}[theorem]{{\bf Corollary}}
\newtheorem{remark}[theorem]{{\bf Remark}}
\newtheorem{conj}[theorem]{{\bf Conjecture}}
\newcommand{\pf}{\smallskip\noindent {\em Proof.}\ \  }
\newcommand{\qed}{\hfill $\Box$\medskip}

\newcommand{\inv}{^{-1}}

\begin{abstract}
Let ${\cal I}_n$ be the symmetric inverse semigroup on $X_n = \{1,
2, \cdots , n\}$ and let ${\cal DP}_n$ and ${\cal ODP}_n$ be its
subsemigroups of partial isometries and of order-preserving partial
isometries of $X_n$, respectively. In this paper we investigate the
cycle structure of a partial isometry and characterize the Green's
relations on ${\cal DP}_n$ and ${\cal ODP}_n$. We show that ${\cal
ODP}_n$ is a $0-E-unitary$ inverse semigroup. We also investigate
the cardinalities of some equivalences on ${\cal DP}_n$ and ${\cal
ODP}_n$ which lead naturally to obtaining the order of the
semigroups. \footnote{This work was carried out when the first
named author was visiting Sultan Qaboos University for a 3-month
research visit in Fall 2010.}\end{abstract}


\section{Introduction and Preliminaries}

Let $X_n=\{1,2, \cdots, n\}$ and ${\cal I}_n$ be the partial
one-to-one transformation semigroup on $X_n$ under composition of
mappings. Then ${\cal I}_n$ is an {\em inverse} semigroup (that is,
for all $\alpha \in {\cal I}_n$ there exists a unique $\alpha' \in
{\cal I}_n$ such that $\alpha = \alpha\alpha'\alpha$ and $\alpha' =
\alpha'\alpha\alpha'$). The importance of ${\cal I}_n$ (more
commonly known as the symmetric inverse semigroup or monoid) to
inverse semigroup theory may be likened to that of the symmetric
group ${\cal S}_n$ to group theory. Every finite inverse semigroup
$S$ is embeddable in ${\cal I}_n$, the analogue of Cayley's theorem
for finite groups, and to the regular representation of finite
semigroups. Thus, just as the study of symmetric, alternating and
dihedral groups has made a significant contribution to group theory,
so has the study of various subsemigroups of ${\cal I}_n$, see for
example \cite{Bor, Fer1, Fer2, Gar, Uma1}.

\noindent A transformation $\alpha \in {\cal I}_n$ is said to be
{\em order-preserving (order-reversing)} if $(\forall x,y \in \
{Dom\, \alpha})\ x \leq y \implies x\alpha \leq y\alpha \ (x\alpha
\geq y\alpha)$ and, is said to be an {\em isometry (or
distance-preserving)} if ($\forall x,y \in {Dom\, \alpha})\, \mid
x-y\mid = \mid x\alpha -y\alpha\mid$. Semigroups of partial
isometries on more restrictive but richer mathematical structures
have been studied \cite{Bra, Wal}. This paper investigates the
algebraic and combinatorial properties of ${\cal DP}_n$ and ${\cal
ODP}_n$, the semigroups of partial isometries and of partial
order-preserving isometries of an $n-$chain, respectively.

In this section we introduce basic terminologies and some
preliminary results concerning the cycle structure of a partial
isometry of $X_n$. In the next section, (Section 2) we characterize
the classical Green's relations and show that ${\cal ODP}_n$ is a
0-E-unitary inverse semigroup. We also show that certain Rees factor
semigroups of ${\cal ODP}_n$ are 0-E-unitary and categorical inverse
semigroups. In Section 3 we obtain the cardinalities of two
equivalences defined on ${\cal DP}_n$ and ${\cal ODP}_n$. These
equivalences lead to formulae for the order of ${\cal DP}_n$ and
${\cal ODP}_n$ as well as new triangles of numbers not yet recorded
in \cite{Slo}.

For standard concepts in semigroup and symmetric inverse semigroup
theory, see for example \cite{How2, Lim, Law1}. Let
\begin{eqnarray} \label{eqn1.1} {\cal DP}_n= \{\alpha \in {\cal I}_n:
(\forall x, y \in X_n)  \mid x-y\mid = \mid x\alpha-y\alpha\mid \}
\end{eqnarray}
\noindent be the subsemigroup of ${\cal I}_n$ consisting of all
partial isometries of $X_n$. Also let
\begin{eqnarray} \label{eqn1.2} {\cal ODP}_n= \{\alpha \in {\cal DP}_n:
(\forall x, y \in X_n)\, x\leq y \Longrightarrow x\alpha\leq
y\alpha\}
\end{eqnarray}
\noindent be the subsemigroup of ${\cal DP}_n$ consisting of all
order-preserving partial isometries of $X_n$. It is clear that if
$\alpha \in {\cal DP}_n$ ($\alpha \in {\cal ODP}_n$) then
$\alpha^{-1} \in {\cal DP}_n$ ($\alpha^{-1} \in {\cal ODP}_n$) also.
Thus we have the following result.

\begin{lemma} \label{lem1.1} ${\cal DP}_n$ and ${\cal ODP}_n$ are
inverse subsemigroups of ${\cal I}_n$.
\end{lemma}

Next we prove a sequence of lemmas that help us understand the cycle
structure of partial isometries. These lemmas also seem to be useful
in investigating the combinatorial questions in Section 3. First,
let $\alpha$ be in ${\cal I}_n$. Then the {\em height} of $\alpha$
is $h(\alpha)= \mid{Im\,\alpha}\mid$, the {\em right [left] waist}
of $\alpha$ is $w^+(\alpha) = max({Im \, \alpha})\, [w^-(\alpha) =
min({Im \,\alpha})]$, the {\em right [left] shoulder} of $\alpha$ is
$\varpi^+(\alpha) = max({Dom \, \alpha})$\, [$\varpi^-(\alpha) =
min({Dom \, \alpha})]$, and {\em fix} of $\alpha$ is denoted by
$f(\alpha)$, and defined by $f(\alpha)=|F(\alpha)|$, where
$$F(\alpha) = \{x \in X_n: x\alpha = x\}.$$

\def\pf{\noindent {\em Proof.}\ \ }

\begin{lemma} \label{lem1.2} Let $\alpha\in {\cal DP}_n$ be such that $h(\alpha)=p$.
Then $f(\alpha)=0 \,or\, 1\, or\, p$.
\end{lemma}

\pf Suppose $x,y \in F(\alpha)$. Then $x=x\alpha$ and $y=y\alpha$.
Let $z\in {Dom\, \alpha}$ where we may without loss of generality
assume that $x<y<z$. Essentially, we consider two cases: $y<z\alpha$
and $x<z\alpha<y$. In the former, we see that $$z-y=\mid z\alpha -
y\alpha\mid=\mid z\alpha - y\mid=z\alpha-y \Longrightarrow
z=z\alpha.$$ \noindent In the latter, we see that $$z-x=\mid z\alpha
- x\alpha\mid=\mid z\alpha - x\mid=z\alpha-x \Longrightarrow
z=z\alpha.$$ However, note that $$\alpha=\pmatrix{2&3&\dots&p+1\cr
1&2&\dots &p}\,\,
\mbox{and}\,\,\beta=\pmatrix{\dots&i-1&i&i+1&\dots\cr
\dots&i+1&i&i-1&\dots}$$ \noindent are nonidempotent partial
isometries with $f(\alpha)=0$ and $f(\beta)=1$. \qed

\begin{corollary} \label{cor1} Let $\alpha\in {\cal DP}_n$. If $f(\alpha)=p>1$ then
$f(\alpha)=h(\alpha)$. Equivalently, if $f(\alpha)>1$ then $\alpha$
is an idempotent.
\end{corollary}

\begin{lemma} \label{lem1.3} Let $\alpha\in {\cal DP}_n$. If $1\in F(\alpha)$ or $n\in
F(\alpha)$ then for all $x\in {Dom \alpha}$, we have $x\alpha =x$.
Equivalently, if $1\in F(\alpha)$ or $n\in F(\alpha)$ then $\alpha$
is a partial identity.
\end{lemma}

\pf Suppose $1 \in F(\alpha)$. Then for all $x\in {Dom \alpha}$,
$x-1=x\alpha-1\alpha=x\alpha-1 \Longrightarrow x=x\alpha$.
Similarly, if $n \in F(\alpha)$, then for all $x\in {Dom \alpha}$,
$n-x=n\alpha-x\alpha=n-x\alpha \Longrightarrow x=x\alpha$. \qed

\begin{lemma} \label{lem1.5} Let $\alpha\in {\cal ODP}_n$ and $n\in {Dom\, \alpha}\cap
{Im\, \alpha}$. Then $n\alpha =n$.
\end{lemma}

\pf Since $n=max({Dom\, \alpha})$ and $n=max({Im\, \alpha})$, and
$\alpha$ is order-preserving then $n\alpha=n$. However, note that in
${\cal DP}_n$ we have $\alpha=\pmatrix{1&n\cr n&1}$, where $n\in
{Dom\, \alpha}\cap {Im\, \alpha}$ but $n\alpha\neq n$.\qed

\begin{lemma} \label{lem1.6} Let $\alpha\in {\cal ODP}_n$ and $f(\alpha)\geq 1$.
Then $\alpha$ is an idempotent.
\end{lemma}

\pf Let $x$ be a fixed point of $\alpha$ and suppose $y\in {Dom
\alpha}$. If $x < y$ then by the order-preserving and isometry
properties we see that $y-x=y\alpha-x\alpha=y\alpha-x
\Longrightarrow y=y\alpha$. The case $y<x$ is similar. However, note
that in ${\cal DP}_n$ we have $\alpha=\pmatrix{1&2\cr 3&2}$, where
$f(\alpha)=1$ but $\alpha^2\neq \alpha$.\qed

\section{Green's relations}

For the definitions of Green's relations we refer the reader to
Howie \cite[Chapter 2]{How2}. It is now customary that when one
encounters a new class of semigroups, one of the questions that is
often asked concerns the characterization of Green's relations. By
Lemma \ref{lem1.1} and \cite[Proposition 2.4.2 \& Ex. 5.11.2]{How2}
we deduce the following lemma.

\begin{lemma} Let $\alpha ,\beta \in {\cal DP}_n.$
Then
\begin{itemize}
\item[(1)] $\alpha \leq _{{\cal R}}\beta $ \,  if and only if\,
${Dom\,\alpha}\subseteq {Dom\,\beta}$;
\item[(2)] $\alpha \leq _{{\cal L}}\beta $\,   if and only if  $Im\,\alpha
\subseteq Im\,\beta $;
\item[(3)] $\alpha \leq _{{\cal H}}\beta $ \,  if and only if\,
${Dom\,\alpha}\subseteq {Dom\,\beta}$ and $Im\,\alpha \subseteq
Im\,\beta $.
\end{itemize}
\end{lemma}

\begin{theorem} Let $S = {\cal DP}_n$ be as defined in (\ref{eqn1.1}). Then
$\alpha \leq _{{\cal D}}\beta $ if and only if there exists an
isometry $\theta : {Dom\,\alpha}\rightarrow {Dom\,\beta}.$
\end{theorem}

\pf \thinspace\ Let $\alpha \leq _{{\cal D}%
}\beta ,$ then there exists $\delta $ in ${\cal DP}_n$ such that
$\alpha \leq _{\cal R}\delta \Longleftrightarrow
\alpha=\delta\eta_1$ \ and $\delta \leq _{\cal L}\beta
\Longleftrightarrow \delta=\eta_2\beta$. Thus
$\alpha=\delta\eta_1=\eta_2\beta\eta_1$ and so ${Dom\,\alpha}
\subseteq {Dom\,\eta_2}$. It is clear that $\eta_2 |_{Dom\,\alpha}$
is an isometry from ${Dom\,\alpha}$ into ${Dom\,\beta}$.

\hspace*{2mm} Conversely, suppose $\theta$ is an isometry from
${Dom\,\alpha}$ into ${Dom\,\beta}$. Define $\eta_1$ by
$x\eta_1=x\theta^{-1}\alpha\,(x\in {Dom\,\beta})$. Then $\eta_1\in
{\cal DP}_n$ and $\theta\eta_1=\theta(\theta^{-1}\alpha)=\alpha$.
Hence $\alpha \leq _{{\cal R}}\theta$. Similarly, define $\eta_2$ by
$x\eta_2=x\theta\beta^{-1}\,(x\in {Dom\,\alpha})$. Then $\eta_2\in
{\cal DP}_n$ and $\eta_2\beta=(\theta\beta^{-1})\beta=\theta$. Hence
$\theta \leq _{{\cal L}}\beta$, as required. \qed

\noindent The corresponding result for ${\cal ODP}_n$ can be proved
similarly.

\begin{theorem} Let $S = {\cal ODP}_n$ be as defined in (\ref{eqn1.2}).
Then $\alpha \leq _{{\cal D}}\beta $ if and only if there exists an
order-preserving isometry $\theta : {Dom\,\alpha}\rightarrow
{Dom\,\beta}.$
\end{theorem}

\noindent Let $E'=E \setminus {0}$. A semigroup S is said to be
$0-E-unitary$ if $(\forall e\in E')(\forall s\in S)\,\,\, es\in E'
\Longrightarrow s\in E'$. The structure theorem for 0-E-unitary
inverse semigroup was given by Lawson \cite{Law2}, see also Szendrei
\cite{Sze} and Gomes and Howie \cite{Gom}. The next result came as a
pleasant surprise to us in the sense that we get a natural class of
0-E-unitary inverse semigroups.

\begin{theorem} \label{thrm1} ${\cal ODP}_n$ is a $0-E-unitary$ inverse
subsemigroup of ${\cal I}_n$.
\end{theorem}

\pf Let $\epsilon\in E({\cal ODP}_n)\setminus \{0\}$. Then  ${Dom\,
\epsilon}\neq \emptyset$ and for $\beta\in {\cal ODP}_n$ such that
$\epsilon\beta$ is a nonzero idempotent we see that ${Dom\,
\epsilon\beta}=\{x, \cdots \}\neq \emptyset$. Now, since
${Dom\,\epsilon\beta} \subseteq {Dom\,\epsilon}$ it follows that
$x\epsilon=x \Longrightarrow x\beta = x\epsilon\beta=x$. Thus, for
any $y\in {Dom \beta}$, (i) if $x<y$, we see that
$y-x=y\beta-x\beta=y\beta-x$. Hence $y=y\beta$, showing that $\beta$
is idempotent, as required. The case (ii) $y<x$ is similar.\qed

\begin{remark} Note that ${\cal DP}_n$ is not 0-E-unitary:
$$\pmatrix{1&2\cr 1&2}\pmatrix{2&3\cr
2&1}=\pmatrix{2\cr 2}\in E({\cal
DP}_n)\,\,\mbox{but}\,\,\pmatrix{2&3\cr 2&1}\notin E({\cal DP}_n).$$
\end{remark}

\noindent For natural numbers $n,p$ with $n\geq p\geq 0$, let
\begin{eqnarray}  L(n,p)=\{\alpha\in {\cal ODP}_n:h(\alpha)\leq p\}
\end{eqnarray}
\noindent be a two-sided ideal of ${\cal ODP}_n$, and for $p>0$, let
\begin{eqnarray} \label{eqn2.1} Q(n,p)= L(n,p)/L(n,p-1)
\end{eqnarray}
\noindent be its Rees quotient semigroup. Then $Q(n,p)$ is a
0-E-unitary inverse semigroup whose nonzero elements may be thought
of as the elements of ${\cal ODP}_n$ of height $p$. The product of
two elements of $Q(n,p)$ is 0 whenever their product in ${\cal
ODP}_n$ is of height less than $p$.

\noindent A semigroup S is said to be $categorical$ if $$(\forall
a,b,c\in S),\,\, abc=0 \Longrightarrow ab=0\,\mbox{or}\, bc=0.$$
\noindent The structure theorem for 0-E-unitary categorical inverse
semigroup was given by Gomes and Howie \cite{Gom}. Now we have

\begin{theorem} \label{thrm2} Let $Q(n,p)$ be as defined in (\ref{eqn2.1}). Then
$Q(n,p)$ is a $0-E-unitary$ categorical inverse semigroup.
\end{theorem}

\pf Let $\alpha, \beta$ and $\gamma\in Q(n,p)$. Note that it
suffices to prove that if $\alpha\beta\neq 0$ and $\beta\gamma\neq
0$ then $\alpha\beta\gamma\neq 0$. Now suppose $\alpha\beta\neq 0$
and $\beta\gamma\neq 0$. Then
${Im\,\alpha\beta}={Im\,\beta}={Dom\,\gamma}$. Hence
$\alpha\beta\gamma\neq 0$, as required.\qed

\begin{remark} Note that ${\cal ODP}_n$ is not categorical: $$\pmatrix{1&2\cr
1&2}\pmatrix{2&3\cr 2&3}\pmatrix{1&3\cr 1&3}= 0$$ \noindent but
$$\pmatrix{1&2\cr 1&2}\pmatrix{2&3\cr 2&3}= \pmatrix{2\cr 2}\neq 0
\,\,\,\mbox{and}\,\,\,\pmatrix{2&3\cr 2&3}\pmatrix{1&3\cr 1&3}=
\pmatrix{3\cr 3}\neq 0.$$
\end{remark}

\section{Combinatorial results}

Enumerative problems of an essentially combinatorial nature arise
naturally in the study of semigroups of transformations. Many
numbers and triangle of numbers regarded as combinatorial gems like
the Stirling numbers \cite[pp. 42 \& 96]{How2}, the factorial
\cite{Mun1, Uma1}, the binomial \cite{Gar}, the Fibonacci number
\cite {How1}, Catalan numbers \cite{Gan}, Lah numbers \cite{Gan,
Lar}, etc., have all featured in these enumeration problems. For a
nice survey article concerning combinatorial problems in the
symmetric inverse semigroup and some of its subsemigroups we refer
the reader to Umar \cite{Uma2}. These enumeration problems lead to
many numbers in Sloane's encyclopaedia of integer sequences
\cite{Slo} but there are also others that are not yet or have just
been recorded in \cite{Slo}.

\noindent Now recall the definitions of {\em height} and {\em fix}
of $\alpha\in {\cal I}_n$ from the paragraph after Lemma
\ref{lem1.1}. As in Umar \cite{Uma2}, for natural numbers $n\geq
p\geq m\geq 0$ we define

\begin{eqnarray} \label{eqn3.1} F(n;p)= \mid\{\alpha \in S:
h(\alpha)=\mid {Im\,\alpha}\mid =p \}\mid,
\end{eqnarray}

\begin{eqnarray} \label{eqn3.2} F(n;m)= \mid\{\alpha \in S:
f(\alpha)=m \}\mid
\end{eqnarray}

\noindent where $S$ is any subsemigroup of ${\cal I}_n$. Also, let
$i=a_i=a$, for all $a \in \{p,m\}$, and $0\leq i\leq n$.

\begin{lemma}\label{lem3.1} Let $S={\cal ODP}_n$. Then
$F(n;p_1)=F(n;1)=n^2$ and $F(n;p_n)=F(n;n)=1$, for all $n\geq 2$.
\end{lemma}

\pf Since all partial injections of height 1 are vacuously partial
isometries, the first statement of the lemma follows immediately.
For the second statement, it is not difficult to see that there is
exactly one partial isometry of height $n$: $\pmatrix{1&2&\dots&n\cr
1&2&\dots &n}$ (the identity). \qed

\begin{lemma} \label{lem3.2} Let $S={\cal ODP}_n$. Then
$F(n;p_2)=F(n;2)=\frac{1}{6}n(n-1)(2n-1)$, for all $n\geq 2$.
\end{lemma}

\pf First, we say that 2-subsets of $X_n$ (that is, subsets of size
2) say, $A=\{a_1,a_2\}$ and $B=\{b_1,b_2\}$ are of the same {\em
type} if $\mid a_1-a_2\mid = \mid b_1-b_2\mid$. Now observe that if
$\mid a_1-a_2\mid = i\,\,(1\leq i\leq n-1)$ then there are $n-i$
subsets of this type. However, for partial order-preserving
isometries once we choose a 2-subset as a domain then the possible
image sets must be of the same type and there is only one possible
order-preserving bijection between any two 2-subsets of the same
type. It is now clear that
$F(n;p_2)=F(n;2)=\sum_{i=1}^{n-1}(n-i)^2=\frac{1}{6}n(n-1)(2n-1)$,
as required. \qed

\begin{lemma} \label{lem3.3} Let $S={\cal ODP}_n$. Then
$F(n;p)=F(n-1;p-1)+F(n-1;p)$, for all $n\geq p\geq 3$.
\end{lemma}

\pf Let $\alpha \in {\cal ODP}_n$ and $h(\alpha)=p$. Then it is
clear that $F(n;p)=\\ \mid A\mid + \mid B\mid$, where $A=\{\alpha
\in {\cal ODP}_n: h(\alpha)=p\,\,\, \mbox{and} \,\, n\notin {Dom\,
\alpha}\cup {Im\, \alpha}\}$ and $B=\{\alpha \in {\cal ODP}_n:
h(\alpha)=p\,\,\, \mbox{and}\,\, n\in {Dom\, \alpha}\cup {Im\,
\alpha}\}$. Define a map $\theta : \{\alpha \in {\cal ODP}_{n-1}:
h(\alpha)=p\} \rightarrow A$ by $(\alpha)\theta=\alpha'$ where
$x\alpha'=x\alpha\,(x\in {Dom\, \alpha}$. This is clearly a
bijection since $n\notin {Dom\, \alpha}\cup {Im\, \alpha}$. Next,
recall the definitions of  $\varpi^+ (\alpha)$ and $w^+(\alpha)$
from the paragraph after Lemma \ref{lem1.1}. Now, define a map $\Phi
: \{\alpha \in {\cal ODP}_{n-1}: h(\alpha)=p-1\} \rightarrow B$ by
$(\alpha)\Phi=\alpha'$ where

\noindent {\bf (i)} $x\alpha' =x\alpha\, (x\in {Dom\,\alpha})\,
\mbox{and}\, n\alpha' =n\,$ (if $\varpi^+ (\alpha)$ = $w^+(\alpha)$
);

\noindent {\bf (ii)} $x\alpha' =x\alpha\, (x\in {Dom\,\alpha})\,
\mbox{and}\, n\alpha' =n-\varpi^+ (\alpha)+w^+(\alpha)<n$\, (if
$\varpi^+ (\alpha)$ \,$>$ $w^+(\alpha)$);

\noindent {\bf (iii)} $x(\alpha')^{-1} =x\alpha^{-1}\, (x\in
{Im\,\alpha})\, \mbox{and}\, n(\alpha')^{-1} =n-\varpi^+
(\alpha)^{-1}+w^+(\alpha^{-1})<n$\, (if $\varpi^+ (\alpha)$ \,$<$
$w^+(\alpha)$).

\noindent In all cases $h(\alpha')=p$, and case (i) coincides with
$n\in {Dom\,\alpha'} \cap {Im\,\alpha'}$; case (ii) coincides with
$n\in {Dom\,\alpha'} \setminus {Im\,\alpha'}$; case (iii) coincides
with $n\in {Im\,\alpha'} \setminus {Dom\,\alpha'}$. Thus $\Phi$ is
onto. Moreover, it is not difficult to see that $\Phi$ is
one-to-one. Hence $\Phi$ is a bijection, as required. This
establishes the statement of the lemma. \qed

\begin{prop}\label{prop1} Let $S = {\cal ODP}_n$ and $F(n;p)$
be as defined in (\ref{eqn1.2}) and (\ref{eqn3.1}), respectively.
Then $F(n;p)= \frac{(2n-p+1)}{p+1}{n\choose p}$, where $n\geq p \geq
2$.
\end{prop}

\pf (The proof is by induction). \newline Basis step: First, note
that $F(n;1)$, $F(n;n)$ and $F(n;2)$ are true by Lemmas \ref{lem3.1}
and \ref{lem3.2}.
\newline Inductive step: Suppose $F(n-1;p)$ is true for all $n-1\geq p$.
(This is the induction hypothesis.) Now using Lemma \ref{lem3.3}, we
see that
\begin{eqnarray*} F(n;p)& = & F(n-1;p-1)+F(n-1;p)\\
& = & \frac{(2n-p)}{p}{n-1\choose
p-1}+\frac{(2n-p-1)}{p+1}{n-1\choose p}
\,\,(by\,\, ind. \,\, hyp.)\\
& = & \frac{(2n-p)}{p}\frac{p}{n}{n\choose
p}+\frac{(2n-p-1)}{p+1}\frac{(n-p)}{n}{n\choose p}\\
& = & \frac{(2n-p)(p+1)+(2n-p-1)(n-p)}{n(p+1)}{n\choose p}\\
& = & \frac{(2n^2-np+n)}{n(p+1)}{n\choose p} =
\frac{(2n-p+1)}{p+1}{n\choose p},
\end{eqnarray*} as required. \qed

\begin{lemma} \label{lem3.4} For integers $n,\,p$ such that
$n\geq p\geq 2$, we have $\sum_{p=2}^{n}\frac{2n-p+1}{p+1}{n\choose
p}= 3\cdot2^n-n^2-2n-3.$
\end{lemma}

\pf It is enough to observe that $2n-p+1=(2n-2p)+(p+1)$. \qed

\begin{theorem}\label{thrm3.1} Let ${\cal ODP}_n$ be as defined in (\ref{eqn1.2}).
Then $$\mid {\cal ODP}_n\mid = 3\cdot2^n-2(n+1).$$
\end{theorem}

\pf It follows from Proposition \ref{prop1} and Lemma \ref{lem3.4},
and some algebraic manipulation. \qed

\begin{lemma} \label{lem3.5} Let $S={\cal ODP}_n$. Then
$F(n;m)={n\choose m}$, for all $n\geq m\geq 1$.
\end{lemma}

\pf It follows directly from Lemma \ref{lem1.6}. \qed

\begin{prop}\label{prop2} Let $S = {\cal ODP}_n$ and $F(n;m)$
be as defined in (\ref{eqn1.2}) and (\ref{eqn3.2}), respectively.
Then $F(n;0)= 2^{n+1}-(2n+1)$.
\end{prop}

\pf It follows from Theorem \ref{thrm3.1}, Lemma \ref{lem3.5} and
the fact that $\mid {\cal ODP}_n\mid=\sum_{m=0}^{n}F(n;m)$. \qed

\begin{remark} The triangles of numbers $F(n;p)$ and $F(n;m)$, the sequence
$F(n;m_0)$ are as at the time of submitting this paper not in Sloane
\cite{Slo}. However, $\mid {\cal ODP}_n\mid$ is \cite[A097813]{Slo}.
For some computed values of $F(n;p)$ and $F(n;m)$ in ${\cal ODP}_n$,
see Tables 3.1 and 3.2.
\end{remark}

\begin{center}

$$\begin{array}{|c|c|c|c|c|c|c|c|c|c|}
\hline
 \,\,\,\,\,n{\backslash}p&0&1&2&3&4&5&6&7&\sum F(n;p)=\mid {\cal ODP}_n \mid
\\ \hline 0&1&&&&&&&&1
 \\ \hline 1&1&1&&&&&&&2
 \\ \hline 2&1&4&1&&&&&&6
\\
\hline 3&1&9&5&1&&&&&16
\\
\hline 4&1&16&14&6&1&&&&38
\\
\hline 5&1&25&30&20&7&1&&&84
\\
\hline 6&1&36&55&50&27&8&1&&178
\\
\hline 7&1&49&91&105&77&35&9&1&368
\\
\hline
\end{array}$$
\end{center}

\begin{center}
Table 3.1
\end{center}

\begin{center}

$$\begin{array}{|c|c|c|c|c|c|c|c|c|c|}
\hline
 \,\,\,\,\,n{\backslash}m&0&1&2&3&4&5&6&7&\sum F(n;m)=\mid {\cal ODP}_n \mid
\\ \hline 0&1&&&&&&&&1
 \\ \hline 1&1&1&&&&&&&2
 \\ \hline 2&3&2&1&&&&&&6
\\
\hline 3&9&3&3&1&&&&&16
\\
\hline 4&23&4&6&4&1&&&&38
\\
\hline 5&53&5&10&10&5&1&&&84
\\
\hline 6&115&6&15&20&15&6&1&&178
\\
\hline 7&241&7&21&35&35&21&7&1&368
\\
\hline
\end{array}$$
\end{center}

\begin{center}
Table 3.2
\end{center}

\begin{remark} \label{rem3} For $p=0,1$ the concepts of
order-preserving and order-reversing coincide but distinct
otherwise. However, there is a bijection between the two sets for
$p\geq 2$, see \cite[page 2, last paragraph]{Fer2}.
\end{remark}

\begin{lemma} \label{lem3.6} Let $\alpha\in {\cal DP}_n$. Then
$\alpha$ is either order-preserving or order-reversing.
\end{lemma}

\pf If $h(\alpha)=2$ then the result is obvious. However, if
$h(\alpha)>2$ we must consider cases. First suppose that
$\{a_1,a_2,a_3\} \subseteq {Dom \alpha}$, where $a_i\alpha=b_i
\,(i=1,2,3)$ and $1\leq a_1<a_2<a_3\leq n$. There are four cases to
consider if $\alpha$ is neither order-preserving or order-reversing:
$b_1<b_3<b_2$, $b_2<b_1<b_3$, $b_2<b_3<b_1$ and $b_3<b_1<b_2$. In
the first case, note that $b_2-b_1=(b_2-b_3)+(b_3-b_1)$. But
$a_3-a_1=(a_3-a_2)+(a_2-a_1)= \mid a_3-a_2\mid+\mid a_2-a_1 \mid \\=
\mid b_3-b_2\mid+\mid b_2-b_1\mid= \mid b_3-b_2\mid+\mid
b_2-b_3\mid+ \mid b_3-b_1\mid=2\mid b_3-b_2\mid+\mid
b_3\alpha^{-1}-b_1\alpha^{-1}\mid= 2\mid b_3-b_2\mid+\mid
a_3-a_1\mid=2\mid b_3-b_2\mid+ a_3-a_1$, which implies that $\mid
b_3-b_2\mid=0 \Leftrightarrow b_3=b_2$. This is a contradiction. The
other three cases are similar. \qed

\noindent We now use Remark \ref{rem3} and Lemma \ref{lem3.6} to
deduce corresponding results for ${\cal DP}_n$ from those of ${\cal
ODP}_n$ above.

\begin{lemma}\label{lem3.7} Let $S={\cal DP}_n$. Then
$F(n;p_1)=F(n;1)=n^2$ and $F(n;p_n)=F(n;n)=2$, for all $n\geq 2$.
\end{lemma}

\begin{lemma} \label{lem3.8} Let $S={\cal DP}_n$. Then
$F(n;p_2)=F(n;2)=\frac{1}{3}n(n-1)(2n-1)$, for all $n\geq 2$.
\end{lemma}

\begin{lemma} \label{lem3.9} Let $S={\cal DP}_n$. Then
$F(n;p)=F(n-1;p-1)+F(n-1;p)$, for all $n\geq p\geq 3$.
\end{lemma}

\begin{prop}\label{prop3} Let $S = {\cal DP}_n$ and $F(n;p)$
be as defined in (\ref{eqn1.1}) and (\ref{eqn3.1}), respectively.
Then $F(n;p)= \frac{2(2n-p+1)}{p+1}{n\choose p}$, where $n\geq p
\geq 2$.
\end{prop}

\begin{theorem}\label{thrm3.2} Let ${\cal DP}_n$ be as defined in (\ref{eqn1.1}).
Then $$\mid {\cal DP}_n\mid = 3\cdot2^{n+1}-(n+2)^2-1.$$
\end{theorem}

\pf It follows from Proposition \ref{prop3}, Lemma \ref{lem3.4} and
some algebraic manipulation. \qed

\begin{lemma} \label{lem3.10} Let $\alpha\in {\cal DP}_n$. For $1<i<n$,
if $F(\alpha)=\{i\}$ then for all $x\in {Dom\,\alpha}$ we have that
$x+x\alpha=2i$.
\end{lemma}

\pf Let $F(\alpha)=\{i\}$ and suppose $x\in {Dom\,\alpha}$.
Obviously, $i+i\alpha =i+i=2i$. If $x < i$ then $x\alpha >i$, for
otherwise we would have $i-x=\mid i\alpha-x\alpha\mid=\\ \mid
i-x\alpha\mid =i-x\alpha\Longrightarrow x=x\alpha$, which is a
contradiction. Thus, $i-x=\\ \mid i\alpha-x\alpha\mid= \mid
i-x\alpha\mid = \mid x\alpha-i\mid=x\alpha-i\Longrightarrow
x+x\alpha=2i$. The case $x>i$ is similar.\qed

\begin{lemma}\label{lem3.11} Let $S={\cal DP}_n$. Then
$F(n;m)={n\choose m}$, for all $n\geq m\geq 2$.
\end{lemma}

\pf It follows from Corollary \ref{cor1}.\qed

\begin{prop}\label{prop4} Let $S={\cal DP}_n$. Then
$F(2n;m_1)=F(2n;1)=\frac{2(2^{2n}-1)}{3}$ and
$F(2n-1;m_1)=F(2n-1;1)=\frac{2(2^{2n-2}-1)}{3}+2^{2n-2}$, for all
$n\geq 1$.
\end{prop}

\pf Let $F(\alpha)=\{i\}$. Then by Lemma \ref{lem3.10}, for any
$x\in {Dom\,\alpha}$ we have $x+x\alpha=2i$. Thus there $2i-2$
possible elements for ${Dom\,\alpha}: (x,x\alpha)\in \{(1,2i-1),
(2,2i-2), \cdots (2i-1,1)\}$. However, (excluding $(i,i)$) we see
that there are $\sum_{j=0}^{2i-2}{2i-2\choose j}=2^{2i-2}$, possible
partial isometries with $F(\alpha)=\{i\}$, where $2i-1\leq n
\Longleftrightarrow i\leq (n+1)/2$. Moreover, by symmetry we see
that $F(\alpha)=\{i\}$ and $F(\alpha)=\{n-i+1\}$ give rise to equal
number of partial isometries. Note that if $n$ is odd the equation
$i=n-i+1$ has one solution. Hence, if $n=2a-1$ we have
$$2\sum_{i=1}^{a-1}2^{2i-2}+2^{2a-2}=\frac{2(2^{2a-2}-1)}{3}+2^{2a-2}$$
\noindent partial isometries with exactly one fixed point; if $n=2a$
we have
$$2\sum_{i=1}^{a}2^{2i-2}=\frac{2(2^{2a}-1)}{3}$$\noindent
partial isometries with exactly one fixed point.\qed

\begin{prop}\label{prop5} Let $S={\cal DP}_n$. Then
$F(n;m_0)=F(n;0)=\frac{13\cdot2^n-(3n^2+9n+10)}{3}$,\\($n\geq 0$, if
$n$ is even) and
$F(n;m_0)=F(n;0)=\frac{25\cdot2^{n-1}-(3n^2+9n+10)}{3}, (n\geq 1,$
if $n$ is odd).
\end{prop}

\pf It follows from Theorem \ref{thrm3.2}, Lemma \ref{lem3.11},
Proposition \ref{prop4} and the fact that $\mid {\cal
DP}_n\mid=\sum_{m=0}^{n}F(n;m)$.\qed

\begin{remark} The triangles of numbers $F(n;p)$ and  $F(n;m)$ and,
the sequences $\mid {\cal DP}_n\mid$ and $F(n;m_0)$, are as at the
time of submitting this paper not in Sloane \cite{Slo}. However,
$F(n;m_1)$ is \cite[A061547]{Slo}. For some computed values of
$F(n;p)$ and $F(n;m)$ in ${\cal DP}_n$, see Tables 3.3 and 3.4.
\end{remark}

\begin{center}

$$\begin{array}{|c|c|c|c|c|c|c|c|c|c|}
\hline
 \,\,\,\,\,n{\backslash}p&0&1&2&3&4&5&6&7&\sum F(n;p)=\mid {\cal DP}_n \mid
\\ \hline 0&1&&&&&&&&1
 \\ \hline 1&1&1&&&&&&&2
 \\ \hline 2&1&4&2&&&&&&7
\\
\hline 3&1&9&10&2&&&&&22
\\
\hline 4&1&16&28&12&2&&&&59
\\
\hline 5&1&25&60&40&14&2&&&142
\\
\hline 6&1&36&110&100&54&16&2&&319
\\
\hline 7&1&49&182&210&154&70&18&2&686
\\
\hline
\end{array}$$
\end{center}

\begin{center}
Table 3.3
\end{center}

\begin{center}

$$\begin{array}{|c|c|c|c|c|c|c|c|c|c|}
\hline
 \,\,\,\,\,n{\backslash}m&0&1&2&3&4&5&6&7&\sum F(n;m)=\mid {\cal DP}_n \mid
\\ \hline 0&1&&&&&&&&1
 \\ \hline 1&1&1&&&&&&&2
 \\ \hline 2&4&2&1&&&&&&7
\\
\hline 3&12&6&3&1&&&&&22
\\
\hline 4&38&10&6&4&1&&&&59
\\
\hline 5&90&26&10&10&5&1&&&142
\\
\hline 6&220&42&15&20&15&6&1&&319
\\
\hline 7&460&106&21&35&35&21&7&1&686
\\
\hline
\end{array}$$
\end{center}

\begin{center}
Table 3.4
\end{center}

\noindent {\bf Acknowledgements}. The first named author would like
to thank Bowen University, Iwo and Sultan Qaboos University for
their financial support and hospitality, respectively.

\vspace{1cm}

\begin{center}
R.\ Kehinde\\
Department of Mathematics and Statistics\\
Bowen University \\
P. M. B. 284,Iwo, Osun State\\
Nigeria.\\
E-mail:{\tt kennyrot2000@yahoo.com}

\end{center}

\begin{center}
A.\ Umar\\
Department of Mathematics and Statistics\\
Sultan Qaboos University \\
Al-Khod, PC 123 -- OMAN\\
E-mail:{\tt aumarh@squ.edu.om}

\end{center}

\end{document}